\input amstex
\documentstyle{amsppt}
\magnification=1200
\hcorrection{0.3in}
\vcorrection{-0.3in}


\define\op#1{\operatorname{#1}}

\NoRunningHeads
\NoBlackBoxes
\topmatter
\title A Generalized $\pi_2$-Diffeomorphism Finiteness Theorem
\endtitle
\date
\enddate
\author Xiaochun Rong  \footnote{A part of this work was done during the author's 2019 summer visit at Capital Normal University, which was partially supported by Beijing Natural Science Foundation Z19003. \hfill{$\,$}} and Xuchao Yao\endauthor
\address Mathematics Department, Rutgers University
New Brunswick, NJ 08903 USA
\endaddress
\email rong\@math.rutgers.edu
\endemail
\address Mathematics Department, Capital Normal University, Beijing,
P.R.C.
\endaddress
\email qiujuli\@vip.sina.com
\endemail

\abstract The $\pi_2$-diffeomorphism finiteness result (\cite{FR1,2}, \cite{PT}) asserts that the diffeomorphic types of compact $n$-manifolds $M$ with vanishing first and second homotopy groups can be bounded above 
in terms of $n$, and upper bounds on the absolute value of sectional curvature and diameter of $M$.
In this paper, we will generalize this $\pi_2$-diffeomorphism finiteness by removing the 
condition that $\pi_1(M)=0$ and asserting the diffeomorphism finiteness on the 
Riemannian universal cover of $M$. 
\endabstract
\endtopmatter
\document

\head 0. Introduction
\endhead

\vskip4mm

A core issue in Riemannian geometry has been controlled topology of a Riemannian 
manifold $M$ by its geometry, in which bounding the number of possible diffeomorphic 
types of $M$ has been a driving force.  

Let's begin with the classical Cheeger's diffeomorphism finiteness result:

\proclaim{Theorem 0.1} {\rm (\cite{Ch})} Given $n ,d, v>0$, there exists a constant, $C(n,d,v)>0$, such 
that compact $n$-manifolds satisfying 
$$|\op{sec}_M|\le 1, \quad \op{diam}(M)\le d,\quad \op{vol}(M)\ge v,$$
contains at most $C(n,d,v)$ many diffeomorphic types. 
\endproclaim

Note that Theorem 0.1 does not hold if one removes any of the bounds (\cite{Ch}). A key is 
to estimate a lower bound on the injectivity radius of $M$ in terms of these bounds. Indeed, the diffeomorphism 
finiteness holds  when weakens the curvature condition to $\op{Ric}_M\ge -(n-1)$ but strengthens 
the volume condition to that $\op{injrad}(M)\ge \rho>0$ (\cite{An}, \cite{AC}). 

By the Gromov's precompactness, the diffeomorphism finiteness follows if one shows that 
any Gromov-Hausdorff convergent sequence of compact $n$-manifolds, $M_i@>\op{GH}>>X$, 
contains a subsequence of same diffeomorphic type. Under the curvature condition, $\op{Ric}_{M_i}
\ge -(n-1)$,  the limit space $X$ is called a non-collapsed Ricci limit space,
and any tangent cone at $x\in X$ is a metric cone (\cite{CC1}). A point $p\in X$ is regular, if its tangent cone
is unique and is isometric to an Euclidean space. If all points on a non-collapsed Ricci limit space 
$X$ are regular (e.g., $\op{injrad}(M_i)\ge \rho>0$), then $X$ is a manifold and for $i$ large, $M_i$ is diffeomorphic 
to $X$ (\cite{CC1,2}). Under the strong regularity that $|\op{sec}_{M_i}|\le 1$, $X$ is indeed a 
$C^{1,\alpha}$-manifold (\cite{Ch},  \cite{GLP}, cf., \cite{GW}, \cite{Pet}) and thus Theorem 0.1 
follows. With $\op{sec}_{M_i}\ge -1$, 
$X$ is a topological manifold which is homeomorphic to $M_i$ (\cite{Pe}),
and for $n\ne 4$, it is a known result in topology that the homeomorphism finiteness implies 
the diffeomophism finiteness (\cite{GPW1,2}, \cite{KiS}).

A different type of diffeomorphism finiteness result, which primarily concerns collapsed manifolds with 
bounded sectional curvature and diameter, was obtained independently in \cite{FR1,2} and \cite{PT}.

\proclaim{Theorem 0.2} {\rm ($\pi_2$-Diffeomorphism finiteness)} Given $n ,d>0$, there exists a constant, $C(n,d)>0$, 
such that compact $n$-manifolds  satisfying 
$$|\op{sec}_M|\le 1, \quad \op{diam}(M)\le d,\quad \pi_1(M)=\pi_2(M)=0,$$
contains at most $C(n,d)$ many diffeomorphic types. 
\endproclaim

Note that Theorem 0.2 fails, if one drops the restriction that $\pi_2(M)=0$; there are
simply connected $7$-manifolds with bounded positive sectional curvature
and  diameter (\cite{Es}). If one removes the assumption that $\pi_1(M)=0$, then odd-dimensional 
spherical space forms are of infinitely many (with vanishing $\pi_2$). In view of this example, 
a naive question is whether the universal covers of manifolds of vanishing $\pi_2$ are of at most finitely many?

The purpose of this paper is to give a positive answer; which generalizes Theorem 0.2 to the following:

\proclaim{Theorem A} {\rm ($\pi_2$-diffeomorpism finiteness on universal covers)} Given $n ,d>0$, there exists a constant, $C(n,d)>0$, such that compact $n$-manifolds satisfying 
$$|\op{sec}_M|\le 1, \quad \op{diam}(M)\le d,\quad \pi_2(M)=0,$$
its Riemannian universal cover $\tilde M$ contains at most $C(n,d)$ many diffeomorphic types. 
\endproclaim

When restricting to simply connected manifolds, Theorem A implies Theorem 0.2,  while manifolds in 
Theorem A may have non-compact universal covers;  e.g., given any $2$-connected 
$k$-manifold $N$, $3\le k\le n-1$ such that $|\op{sec}_N|\le 1$ and $\op{diam}(N)<d$,  
the metric product with a collapsed flat torus, $N\times \epsilon T^{n-k}$  ($\epsilon>0$ small), 
has a non-compact universal cover. 

In view of Theorems 0.1 and 0.2, the following problem is natural: Is there a diffeomorphism finiteness  of universal covers of 
compact $n$-manifolds satisfying 
$$|\op{sec}_M|\le 1,\quad \op{diam}(M)\le d,\quad \op{vol}(B_1(\tilde p))\ge v>0,\quad \text{at some $\tilde p\in \tilde M$}.$$

We do not have an answer to the above problem. In the Appendix, we will prove a universal cover diffeomorphism (resp. homeomorphism) finiteness for the class of compact $n$-manifolds of non-negative Ricci (resp. sectional) curvature
whose universal covers are not collapsed (see Theorem 3.1). 

The proof of Theorem A adapts the approach to Theorem 0.2 in \cite{FR1,2} (\cite{PT}). To explain 
problems in our proof, let's first briefly describe the proof of Theorem 0.2.

By Theorem 0.1,   it suffices to prove a diffeomorpic stability for a sequence of collapsed manifolds in Theorem 0.2, where the main tool is  the nilpotent fibration structure discovered on collapsed manifolds with bounded sectional curvature (\cite{Gr}, \cite{CG1,2}, \cite{CFG}, \cite{Fu1-3}); the following fibration theorem is a strengthened 
version (see Section 1):

\proclaim{Theorem 0.3} {\rm (Nilpotent bundles, \cite{Gr}, \cite{Fu1,2}, \cite{CFG})} Let a sequence of compact $n$-manifolds, $M_i@>\op{GH}>>X$, such that $|\op{sec}_{M_i}|\le 1$ and $X$ is compact with $\dim(X)<n$. 
Passing to a subsequence, the sequence of the frame bundles, 
$(F(M_i),O(n))$, equipped with canonical metrics, equivariant converges to $(Y,O(n))$, and $Y$ is a $C^{1,\alpha}$-manifold. Then for $i$ large, there is a smooth $O(n)$-invariant bundle map, $f_i: (F(M_i),O(n))\to (Y,O(n))$,  satisfying 
the following properties:   

\noindent {\rm(0.3.1)} $f_i$ is an $\epsilon_i$-GHA i.e. $f_i$ is an $\epsilon_i$-isometry, $\epsilon_i\to 0$.

\noindent {\rm(0.3.2)} A fiber $F_i$ is diffeomorphic to an nilmanifold, and the structural group preserves affine
structures on fibers.

\noindent {\rm(0.3.3)} Norms of the second fundamental form of $f_i$-fibers are uniformly bounded by $c(n)$. 
\endproclaim

The $O(n)$-invariance implies that $f_i$ descends to a possibly singular fibration map,  $\bar f_i: M_i\to X=Y/O(n)$. Note that if $M_i$ is orientable, then Theorem 0.3 applies to the oriented frame bundle, $\op{SF}(M_i)$, which is a component of $F(M_i)$. 

If $\pi_1(M_i)=0$, then the singular nilpotent structure coincides with an 
almost isometric torus $T^k$-action on $M_i$ (\cite{Ro}), $k=n-\dim(X)$, which is descended from the $\op{SO}(n)$-invariant principal $T^k$-bundle on the oriented frame bundle,  $T^k\to (\op{SF}(M_i),\op{SO}(n))\to (Y,\op{SO}(n))$. 

Two principal $G$-bundles, $G\to E_i@>\pi_i>> B$, are weakly equivalent, if there is a diffeomorphism, $f: E_1\to E_2$, and an automorphism, $A: G\to G$, such that $\pi_1=\pi_2\circ f$ and $f(s(x))=A(s)f(x)$ for all $x\in E$ and $s\in G$.

The following is a standard result from principal bundle classification theory (\cite{FR1}).

\proclaim{Lemma 0.4}  Given two principal torus $T^k$-bundles, $T^k\to E_j\to B$, $j=1, 2$, where 
$E_j$ and $B$ are manifolds, if $\pi_1(E_j)=\pi_2(E_j)=0$, 
then the two principal $T^k$-bundles are weakly equivalent. In particular, $E_1$ is diffeomorphic to $E_2$.
\endproclaim

Note that without loss of generality, we may assume $\dim(M)\ge 3$. Because  $\op{SF}(M_i)$ may not be simply connected (e.g., $\pi_1(\op{SF}(S^3))\cong \Bbb Z_2$), either one 
uses a generalization of Lemma 0.4 (\cite{FR2}), or for a simplicity one may replace $\op{SF}(M_i)$  with the associated $\op{Spin}(n)$-bundle over $M_i$, 
$\op{Spin}(M_i)$, where $\op{Spin}(n)\to \op{SO}(n)$ is the universal cover group (Theorem 0.3 trivially 
extends to this set-up). In the rest of the paper, we will use $S(M)$ (resp. $S(n)$) to denote either $\op{SF}(M)$ (resp. $\op{SO}(n)$) when $\pi_1(\op{SO}(n))$ has a trivial image in $\pi_1(\op{SF}(M))$ or otherwise $\op{Spin}(M)$ (resp. $\op{Spin}(n))$.

If, in addition,  $\pi_2(M_i)=0$, then  $\pi_2(S(M_i))=0$ (because $\pi_2(S(n))=0$, the second homotopy of any Lie group is zero), and by Lemma 0.4  all the principal $T^k$-bundles are weakly equivalent, i.e., fixing $i_0$ and for any $i$, there is a bundle equivalence (may not be principal $T^k$-bundle equivalence),  $f_i: S(M_{i_0})\to S(M_i)$, and an automorphism $\phi_i:  T^k\to T^k$, such that $f_i(t(x))=\phi_i(t)f_i(x)$ for all 
$t\in T^k$ and $x\in S(M_{i_0})$. To conclude that $M_{i_0}$ is diffeomorphic to $M_i$,  one can modify $f_i$ to an $S(n)$-conjugate map, via the center of mass technique (\cite{Pa}, \cite{GK}).  

Let's now remove the assumption that $M_i$ is simply connected, and first try to establish a diffeomorphism
stability on $S(\tilde M_i)$, up to a subsequence. In view of the above, we will study the lifting $S(n)$-invariant 
bundle on $S(\tilde M_i)$ from $F_i\to (S(M_i),S(n)_i)\to (Y,S(n))$, which is a $S(n)$-invariant principal nilpotent Lie group $N_i$-bundle, $N_i\to (S(\tilde M_i),S(n)_i)\to (\tilde Y,S(n))$, where $\tilde Y$ is the Riemannian universal cover of $Y$  with lifting 
$S(n)$-action (see Lemma 2.2). 
Unfortunately, it is not possible to have a weakly equivalence for principal nilpotent Lie group bundles (comparing to Lemma 0.4); because contracting to abelian Lie groups, in each dimension $\ge 3$ there are infinitely 
many isomorphism class of nilpotent Lie groups. We observe that if $N_i$ contains no compact subgroup, i.e., $N_i$ is 
diffeomorphic to $\Bbb R^m$, then the principal $N_i$-bundle is trivial, i.e., $S(\tilde M_i)=\tilde Y_i\times N_i$ are all diffeomorphic (see Lemma 2.5).

We will study the $S(n)$-invariant sub-principal abelian bundles,  
$$T^k\to (S(\tilde M_i),S(n)_i)\to (S(\tilde M_i)/T^k,S(n)_i),$$
where $T^k$ denotes the maximal compact subgroup of $N_i$ and the $S(n)_i$-action on $S(\tilde M_i)/T^k$ is descending from the $S(n)_i$-action on $S(\tilde M_i)$. In order to apply Lemma 0.4, we will prove that the diffeomorphic type of based spaces, $S(M_i)/T^k\overset{\op{diffeo}}\to \cong \tilde Y\times \Bbb R^m$,  $m=\dim(N_i)-k$ (see Lemma 2.5), for all $i$.  To be able to 
modify a weakly equivalence to a $S(n)_i$-conjugate map  (so as to conclude the diffeomorphic stability for $\tilde M_i$), 
it requires that passing to a subsequence the descending $S(n)_i$-actions on $\tilde Y\times \Bbb R^m$ are conjugate.  
We will prove this by induction on the length of the nilpotency (passing to a subsequence, we may assume $N_i$ has the same length of nilpotency  such that the quotient of two adjacent normal nilpotent subgroups are independent of $i$), which relies on that we replace the original metric on $M_i$ by a nearby $N_i$-invariant metric of high regularity (\cite{CFG}) (see Lemma 2.7).  
 
\vskip4mm

The rest of the paper is organized as follows: 

\vskip2mm

In Section 1, we will review basic notions and tools that will be used in the proof of Theorems A.

In Section 2, we will prove Theorem A.

In Appendix, we will prove a universal cover diffeomorphism (resp. homeomorphism) finiteness 
for a class of compact manifolds of non-negative Ricci (resp. sectional) curvature.  

\vskip10mm

\head 1.  Preliminaries
\endhead 

\vskip4mm

We will review the notions of equivariant GH-convergence, a form of stability of isometric compact Lie group actions, 
the singular fibration theorem.

\vskip4mm

\subhead a. Equivariant Gromov-Hausdorff convergences
\endsubhead

\vskip4mm

The reference of this subsection is [FY] (cf. [Ro]).

Let $X_i@>GH>>X$ be a convergent sequence of compact length metric spaces, i.e.,
there is a sequence $\epsilon_i\to 0$ and a sequence of maps,
$h_i: X_i\to X$, such that $|d_X(h_i(x_i),h_i(x_i'))-d_{X_i}(x_i,x_i)|
<\epsilon_i$ ($\epsilon_i$-isometry), and for any $x\in X$, there is $x_i
\in X_i$ such that $d_X(h_i(x_i),x)<\epsilon_i$ ($\epsilon_i$-onto), and
$h_i$ is called an $\epsilon_i$-Gromov-Hausdorff approximation, briefly,
$\epsilon_i$-GHA.

Assume that $X_i$ admits a closed group $\Gamma_i$-action by isometries. Then
$(X_i,\Gamma_i)@>\op{eqGH}>>(X,\Gamma)$ means that there are a closed group $\Gamma$ of
isometries on $X$, a sequence $\epsilon_i\to 0$ and
a sequence of triple of $\epsilon_i$-GHA, $(h_i,\phi_i,\psi_i)$, $h_i: X_i\to X$, $\phi_i:
\Gamma_i\to \Gamma$ and $\psi_i: \Gamma\to \Gamma_i$, such that for all $x_i\in X_i, \gamma_i\in \Gamma_i$ and $\gamma\in \Gamma$,
$$d_X(h_i(x_i), \phi_i(\gamma_i)h_i(\gamma_i^{-1}(x_i)))<\epsilon_i,\quad
d_X(h_i(x_i),\gamma^{-1}(h_i(\psi_i(\gamma)(x_i)))<\epsilon_i,$$
where $\Gamma_i$ and
$\Gamma$ are equipped with the induced metrics from $X_i$ and $X$.
We call $(h_i,\phi_i,\psi_i)$
an $\epsilon_i$-equivariant GHA.

When $X$ is not compact, the above notion of equivariant convergence
naturally extends to a pointed version $(h_i,\phi_i,\psi_i)$: $h_i:
B_{\epsilon_i^{-1}}(p_i)\to B_{\epsilon_i^{-1}+\epsilon_i}(p)$, $h_i(p_i)=p$,
$\phi_i: \Gamma_i(\epsilon_i^{-1})\to \Gamma(\epsilon_i^{-1}+\epsilon_i)$,
$\phi_i(e_i)=e$, $\psi_i: \Gamma(\epsilon_i^{-1})\to \Gamma_i(\epsilon_i^{-1}
+\epsilon_i)$, $\psi_i(e)=e_i$, and the above inequalities hold whenever the multiplications
stay in the domain of $h_i$, where $\Gamma_i(R)=\{\gamma_i\in\Gamma_i,\,
\,d_{X_i}(p_i,\gamma_i(p_i))\le R\}$.

\proclaim{Lemma 1.1}  Let $(X_i,p_i)@>\op{GH}>>(X,p)$, where $X_i$ is a  complete
locally compact length space. Assume that $\Gamma_i$ is a closed group of isometries
on $X_i$. Then there is a closed group $G$ of isometries on $X$
such that passing to a subsequence, $(X_i,p_i,\Gamma_i)@>\op{eqGH}>>(X,p,G)$.
\endproclaim

\proclaim{Lemma 1.2}  Let $(X_i, p_i, \Gamma_i)@>\op{eqGH}>>(X,p,G)$, where $X_i$ is
a complete locally compact length space and $\Gamma_i$ is a closed subgroup
of isometries. Then $(X_i/\Gamma_i,\bar p_i)@>\op{GH}>>(X/G,\bar p)$.
\endproclaim

\vskip4mm

\subhead b. The center of mass
\endsubhead

\vskip4mm

The stability of compact Lie group $G$-actions in \cite{Pa} asserts that says if two $G$-actions on a 
compact manifold $M$ are $C^1$-close, then the two $G$-actions are conjugate via a diffeomorphism 
close and isotopy to $\op{id}_M$. In \cite{GK}, a geometric criterion for the $C^1$-closeness is described. 
For our purpose, we state the following version.

\proclaim{Theorem 1.3} Let $M$ be a complete manifold which admits two compact Lie group $G$-actions,
$\mu_i: G\times M\to M$, $i=0, 1$. Assume that $M$ admits two invariant metrics, $g_i$, with respect to 
$\mu_i$  respectively such that 
$$|\op{sec}_{g_i}|\le 1,\quad \op{injrad}(g_0)\ge i_0.$$
Then there is a constant $\epsilon(n,i_0)>0$ such that if 
$$d((M,G,\mu_0),(M,G,\mu_1))=\sup_{x\in M}\{d_{g_0}(\mu_0(g,x),\mu_1(g,x))\,\, \, g\in G\}<\epsilon(n,i_0),$$
then there is $G$-conjugate diffeomorphism, $f: (M,G,\mu_0)\to (M,G,\mu_1)$ which is $\epsilon$-isotopy to 
$\op{id}_M$.
\endproclaim

Note that the assumption of Theorem 1.3 on injectivity radius can be replaced by $\op{vol}(B_1(x,g_0))\ge v_0>0$ for
all $x\in M$ (cf. \cite{Ch}). 

\vskip4mm

\subhead c.  Singular fibrations
\endsubhead

\vskip4mm

Consider a collapsing sequence of compact $n$-manifolds, $M_i@>\op{GH}>>X$, such that $|\op{sec}_{M_i}|\le 1$ and $\op{diam}(M_i)\le d$. Let $F(M_i)$ denote the orthogonal frame bundles. The Levi-Civita connection on $M$ determines a `horizontal' distribution on $F(M_i)$, and thus an $O(n)$ bi-invariant metric uniquely determines 
a Riemannian metric on $F(M)$ such that the projection map, $F(M_i)\to M$ is a Riemannian submersion. Passing
to a subsequence, one obtains the following equivariant convergence (see Lemma 1.2),
$$\CD (F(M_i),O(n)_i)@>\op{eqGH}>> (Y,O(n))\\
@VV \op{Proj}_i V @VV \op{Proj} V\\
M_i,@>\op{GH}>>X=Y/O(n),
\endCD$$
where $O(n)_i$ denotes the isometric $O(n)$-action on $F(M_i)$; $\mu_i: O(n)\times F(M_i)\to F(M_i)$.
In \cite{Fu1,2}, it was proved that there is an $O(n)$-invariant fibration, $f_i: (F(M_i),O(n))\to (Y,O(n))$, 
as in Theorem 0.3, without that $Y$ is a $C^{1,\alpha}$-manifold; which follows from the following 
curvature estimate for the canonical metric on $F(M_i)$.

\proclaim{Theorem 1.4} {\rm (\cite{KoS})}  Let $M$ be a complete Riemannian $n$-manifold, and let $F(M)$ 
be the frame bundle equipped with a canonical metric such that $F(M)\to M$ is a Riemannian submersion.
Assume that $|\op{sec}_M|\le 1$. The sectional curvature of the canonical metric on $F(M_i)$ satisfies 
$$|\op{sec}_{F(M_i)}|\le c(n).$$
\endproclaim

The proof of Theorem 0.3 in \cite{Fu2} and \cite{CFG} uses a smoothing technique to find, on $M_i$, a nearby
metric of higher regularity so that its canonical metric on $F(M_i)$ converges to a smooth manifold 
which is diffeomorphic to $Y$. Using Theorem 1.4, one may apply the $\delta$-splitting map technique in 
\cite{CC1} to construct a fiber bundle map from $F(M_i)$ to $Y$ (Lemma 3.7). Note that under the condition 
that $|\op{sec}_{F(M_i)}|\le C$, it is easy to see that  $L^2$-estimate for a $\delta$-splitting map (see Definition 3.8) 
implies a pointwise estimate, and thus $f_i$ is non-degenerate, i.e.,
$f_i$ is a fiber bundle map. This will give a direct construction of a fiber bundle map without relying a 
smoothing technique.

\demo{A verification that $Y$ is a $C^{1,\alpha}$-manifold}

For any $y\in Y$, $x_i'\in F(M_i)$  such that 
$f_i(x_i')=y$, let $S_\rho(x_i')=\exp_{x_i'}B^\perp_\rho(0)$ denote a normal slice of 
$f_i^{-1}(y)$, where $B^\perp_\rho(0)\subset T^\perp_{x_i'}F(M_i)$, the normal subspace to 
$f^{-1}_i(y)$, and $\rho=\rho(n)$ is a constant for all $i$ (this is a consequence of local 
bounded covering geometry property, \cite{CFG}).  Equipped with the induced metric, 
$S_\rho(x_i')@>\op{GH}>>B_\rho(y)$. It is easy to see  that the induced metric on $S_\rho(x_i')$ 
has bounded sectional curvature, and there is a uniform lower bound on the injectivity
radii on $S_{\frac \rho2}(x_i)$. By now the $C^{1,\alpha}$-regularity follows (\cite{An}).
\qed\enddemo

\vskip4mm

\head 2. Proof of Theorem A
\endhead 

\vskip4mm

By the Gromov's precompactness, it suffices to consider a sequence of compact $n$-manifolds 
satisfying the conditions of Theorem A, $M_i@>\op{GH}>>X$, and show that passing to a subsequence, all $\tilde M_i$
are diffeomorphic. Without loss of the generality, we may assume that $M_i$ is oritentable (so we will 
use $S(M_i)$ as previous defined).

By Theorem 0.1, we may assume that $\dim(X)<n$.
By Theorem 0.3, we obtain the following equivariant commutative diagram:
$$\CD F_i@> >>(S(M_i),S(n)_i)@> f_i>>(Y,S(n))\\
@VV \op{proj}_iV @VV \op{proj}_i V @VV \op{proj} V\\
\bar F_i@> >> M_i@>\bar f_i>> X=Y/S(n),
\endCD \tag 2.1$$ 
where $f_i$ is a $S(n)$-invariant fiber bundle map and an $\epsilon_i$-GHA ($\epsilon_i\to 0$), $F_i$ is a 
nilpotent manifold, and $\bar f_i$ is a 
singular fibration. 

We point it out that at the end of the proof (see the proof Lemma 2.7), we will use a nearby invariant metric on $M_i$ with higher regularity, 
i.e., the induced metric on each $F_i$-fiber is left invariant, which is obtained by a smoothing technique and used
in the construction of $f_i$  (\cite{CFG}).

\vskip4mm

\subhead d. The lifting principal connected nilpotent Lie group bundles
\endsubhead 

\vskip4mm
 
\proclaim{Lemma 2.2} {\rm (Lifting principal nilpotent bundles)} Let $F_i\to (S(M_i), S(n)_i)@>f_i>>(Y, S(n))$
be the $S(n)$-invariant nilpotent fiber bundle satisfying (2.1). 
Then there is a lifting $S(n)$-invariant principal connected nilpotent  Lie group
$N_i$-bundle on the Riemannian universal cover of $S(M_i)$, $N_i\to (S(\tilde M_i),S(n)_i)@> \tilde f_i>> (\tilde Y,S(n))$, 
with the following commutative diagram:
$$\CD N_i@>  >> (S(\tilde M_i), S(n)_i,\tilde p_i)@>\tilde f_i >>(\tilde Y, S(n),\tilde p)\\
@VV \pi_i|_{N_i} V @VV \pi_i V @VV \pi V\\
F_i@> >> (S(M_i), S(n)_i,p_i)@>f_i >>(Y, S(n),p),
\endCD\tag 2.3$$
where $\pi$ and $\pi_i$ are Riemannian universal covering maps. 
\endproclaim

Note that $\tilde f_i$ may not be an $\epsilon_i$-GHA; e.g., $\op{diam}(N_i(\tilde x_i))\ge 1$ (see Example 2.4).

\demo{Proof}  Let $\pi_i: (\tilde M_i,\tilde p_i)\to (M_i,p_i)$ be the Riemannian universal cover. Then $\pi_i$
induces a Riemannian universal cover, $\pi_i: (S(\tilde M_i),\tilde p_i)\to (S(M_i),p)$. 
Let $\tilde F_i$ denote a component of $\pi_i^{-1}(F_i)$.
Let $\Lambda_i=\op{Im}(\pi_1(F_i)\to \Gamma_i)$, which is a normal subgroup of $\Gamma_i=\pi_1(S(M_i))$, 
and let $\bar \Gamma_i=\Gamma_i/\Lambda_i$.  Then
$$\pi_i^{-1}(F_i)=\bigcup_{\bar \gamma_j\in \bar \Gamma_i} \bar \gamma_j(\tilde F_i), \quad \bar \gamma_j(\tilde F_i)\cap \bar \gamma_l(\tilde F_i)=\cases \emptyset & j\ne l\\  \bar \gamma_i(\tilde F_i) & j=l\endcases.$$

We will call a component of $\pi_i^{-1}(F_i)$ a fiber on $S(\tilde M_i)$, which is diffeomorphic to $\tilde F_i$.
Then $S(\tilde M_i)$ is a disjoint union of $\tilde F_i$-fibers. Because $F_i\to S(M_i)\to Y$ is
a fiber bundle, the fibers on $S(\tilde M_i)$ has a local trivialization and thus form a
fiber bundle., $\tilde F_i\to S(\tilde M_i)\to \tilde Y_i$. It is clear that $\tilde F_i$ and $\tilde Y_i$ are covering spaces of 
$F_i$ and $Y$ respectively. Because $\pi_1(S(\tilde M_i))=0$, $\pi_1(\tilde Y_i)=0$ 
and thus $\tilde Y_i=\tilde Y$ is the Riemannian universal cover of $Y$. By the long homotopy exact sequence,
$$\pi_2(S(\tilde M_i))\to \pi_2(\tilde Y_i)\to \pi_1(\tilde F_i)\to 0=\pi_1(S(\tilde M_i)),$$
we see that $\pi_1(\tilde F_i)$ is abelian. Because the universal cover of $\tilde F_i$ is a simply connected nilpotent Lie group whose center contains $\pi_1(\tilde F_i)$, $\tilde F_i$ is isomorphic to a connected nilpotent Lie group $N_i$. 
\qed\enddemo

\example{Example 2.4} Let $(S^3,g_\epsilon)$ denote the Berger's family of collapsing $3$-sphere.
Then $(M,\bar g_\epsilon)=(S^3,g_\epsilon)\times S^1_\epsilon@>\op{GH}>> S^2_{\frac 12}$. Then the 
fiber bundle in Theorem 0.3 is given by the free isometric $T^2$-action on $(S^3,g_\epsilon)\times S^1_\epsilon$, 
the lifting principal bundle is $S^1\times \Bbb R^1\to (S^3\times \Bbb R^1,S(3))\to (S^2_{\frac 12},S(3))$, where 
$(\tilde M, \tilde g_\epsilon)=(S^3,\times \Bbb R^1,g_\epsilon\otimes ds^2)$ satisfies that $\op{diam}_{\bar g_\epsilon}(S^1\times \Bbb R^1)=\infty$ 
and $\op{vol}_{\bar g_\epsilon}(B_1(\tilde p))\to 0$ as $\epsilon\to 0$. 
\endexample

\vskip4mm

\subhead e. Principal bundles of simply connected nilpotent Lie group are trivial
\endsubhead 

\vskip4mm

Recall that fixing a dimension $k\ge 3$, there are infinitely many isomorphism classes of a connected 
nilpotent Lie groups. Hence, in general one cannot expect Lemma 0.4 for a principal nilpotent Lie group bundles.
As seen in the introduction, we will consider the $S(n)$-invariant principal torus $T^k$-bundle, 
$$T^k\to (S(\tilde M_i),S(n)_i)\to (S(\tilde M_i)/T^k,S(n)_i).$$

The following lemma is a key property to the equivariant stability of the principal connected $T^k$-bundles. 

\proclaim{Lemma 2.5}  Let $N\to M\to Y$ be a  smooth principal connected nilpotent Lie group bundle. If $T^k$ is the maximal compact subgroup of $N$, then principal $T^k$-bundle, $T^k\to M\to M/T^k$,  satisfies that $M/T^k$ is diffeomorphic to the 
product $Y\times \Bbb R^m$, $m=\dim(N)-k$.
\endproclaim

\demo{Proof}  Note that $T^k$ is contained in the center of $N$, and we may assume that $k<\dim(N)$. 
Then $\hat N=N/T^k$ is a simply connected nilpotent Lie group, and we shall show that 
$\hat N\to M/T^k\to Y$ is a trivial principal $\hat N$-bundle. 

We may assume that the principal $\hat N$-bundle is the pull back bundle of a map, $f: Y\to B\hat N$, where 
$\pi: E\hat N\to B\hat N$ is an universal principal $\hat N$-bundle, where $E\hat N$ is characterized as a path connected 
topological space on which $\hat N$ acts freely and $E\hat N$ has only trivial homotopy groups (\cite{MT}). 
It is clear that $E\hat N$ can be chosen to be $\hat N$, and thus $B\hat N=\{\op{pt}\}$. By now, the desired result follows.
\qed\enddemo

\vskip4mm

\subhead f. Modifying a bundle isomorphism to a $S(n)$-conjugate bundle isomorphism
\endsubhead

\vskip4mm

Returning to (2.3), passing to a subsequence we may assume that all compact subgroups of $N_i$ is isomorphic 
to $T^k$. By Lemma 2.5,
we obtain a sequence of $S(n)$-invariant principal $T^k$-bundles,
$$T^k@> >> (S(\tilde M_i),S(n)_i)@> \op{proj}_i>> (\tilde Y\times \Bbb R^m,S(n)_i),$$
where the $S(n)_i$-action on $\tilde Y\times \Bbb R^m$ was descended from the $S(n)_i$-action on $S(\tilde M_i)$ 
(note that the descending $S(n)_i$-action on $\tilde Y$ is independent of $i$; see Theorem 0.3). 
Because $\pi_1(S(\tilde M_i))=\pi_2(S(\tilde M_i))=0$, fixing $i_0$ and for any $i\ge i_0$ by Lemma 0.4 
we may assume a principal $T^k$-bundle weakly equivalence, 
$$\CD T^k@= T^k\\
@VV V   @VV V\\
(S(\tilde M_{i_0}),S(n)_{i_0})@>\phi_i>> (S(\tilde M_i),S(n)_i)\\
@VV \op{proj}_{i_0}V @VV \op{proj}_i V\\
(\tilde Y\times \Bbb R^m,S(n)_{i_0})@=(\tilde Y\times \Bbb R^m,S(n)_i).
\endCD$$
Via $\phi_i$, we identify $S(\tilde M_{i_0})$ with $S(\tilde M_i)$, and let $\tilde g_i$ denote the pullback metrics by $\phi_i^*$  on $S(\tilde M):=S(\tilde M_{i_0})$. Then 
$$\CD T^k@= T^k\\
@VV V   @VV V\\
(S(\tilde M),S(n)_i)@>\op{id}_{S(\tilde M)}>> (S(\tilde M),S(n)_j)\\
@VV \op{proj} V @VV \op{proj}V\\
(\tilde Y\times \Bbb R^m,S(n)_i)@> \op{id}_{\tilde Y\times \Bbb R^m}>>(\tilde Y\times \Bbb R^m,S(n)_j).
\endCD\tag 2.6$$

Note that because $\phi_i$ is a weakly equivalence of bundles, in the above $\op{id}_{S(\tilde M)}$ is a weakly equivalence of $T^k$-bundles (which may not be a  principal $T^k$-bundle equivalence).

In order to modified $\op{id}_{S(\tilde M)}$ to a $S(n)$-conjugate map (see Lemma 2.8), we need the following

\proclaim{Lemma 2.7} Let $T^k\to (S(\tilde M),S(n)_i)@> \op{proj}_i>> (\tilde Y\times \Bbb R^m,S(n)_i)$ be 
the $S(n)$-invariant connected principal abelian Lie group bundle. Passing to  a subsequence, 
all $(\tilde Y\times \Bbb R^m,S(n)_i)$ are conjugate by diffeomorphisms that are close and isotropy to 
$\op{id}_{\tilde Y\times \Bbb R^m}$.
\endproclaim

\demo{Proof}  As pointed out at the beginning of this section, without loss of generality we will assume the following: the metric 
$g_i$ on $M_i$ is invariant, i.e., the induced metric on any $N_i$-fiber is left invariant, and $|\nabla ^sg_i|\le c(n,s)$ for 
all $i$ and $s$ (\cite{CFG}). Moreover, passing to a subsequence, we may assume the following normal descending sequence:
$$N_i/T^k=\hat N_i=\hat N_{i,0}\triangleright \hat N_{i,1}\triangleright\cdots \triangleright \hat N_{i,s}\triangleright e,\quad \hat N_{i,j-1}/\hat N_{i,j}=\Bbb R^{k_j} \,\,\text{(independent of $i$).}$$

We shall  proceed the proof by induction on $s$, starting with $s=0$, i.e., $\hat N_i$ is 
isomorphic to abelian group $\Bbb R^m$. By Lemma 2.5, the $S(n)$-invariant principal $\Bbb R^m$-bundle,
$$\hat N_i\to (S(\tilde M)/T^k,S(n)_i)\to (\tilde Y,S(n)),$$
is trivial, $S(\tilde M)/T^k=\tilde Y\times \hat N_i=\tilde Y\times \Bbb R^m$; which gives a natural identification
on all principal $\Bbb R^m$-bundles.  Note that the abelian Lie group $\Bbb R^m$ acts on the $\Bbb R^m$-factor 
by multiplications and on $\tilde Y\times \Bbb R^m$ by isometrices with respect to $\bar g_i$ ($\bar g_i$ is 
in general not a product metric). We observe the following:

\vskip1mm

\noindent (2.7.1) The $\Bbb R^m$-action on $\tilde Y\times \Bbb R^m$ commutes with the $S(n)_i$-action.

\noindent (2.7.2) Fixing $(\tilde p,0)\in \tilde Y\times \Bbb R^m$, let $D_i\ni (\tilde p,0)$ 
denote a fundamental domain of $S(M_i)/T^k$, then $S(n)_i(\tilde p,0)\subset D_i$ (note that 
$\pi_i: (\tilde Y\times \hat N_i,(\tilde p,0))\to (S(M_i)/T^k,p)$, is a universal covering map).

\noindent (2.7.3) Two quotient metrics on $\tilde Y\times \Bbb R^m$, $\bar g_i$ and $\bar g_j$, has 
bounded norm on all derivatives. Hence, passing to a subsequence we may assume $\bar g_i$ and 
$\bar g_j$  are $C^\infty$-close on a compact subset $C$  (e.g., $C$ is the closed ball, $\bar 
B_{4\op{diam}(S(M_{i_0}))}((\tilde p,0),\bar g_{i_0})$) containing both 
fundamental  domains, $D_i$ and $D_j$), so the isometric $S(n)_i$-action (with respect to $\bar g_i$) 
and the isometric $S(n)_j$-action (with respect to $\bar g_j$) are equivariant close on $C$ (Lemma 1.1).

\vskip1mm

By the $S(n)$-invariance, we see that the set,
$$S(\tilde y,r)=\{s^{-1}_js_i(\tilde y,r),\,\,s\in S(n)\}\subset \tilde y\times \Bbb R^m, \quad s_i=\bar \mu_i(s,\cdot),$$
where $\bar \mu_i$ denotes the $S(n)_i$-action on $\tilde Y\times \Bbb R^m$. 
We define the map, 
$$\psi_{i,j}: S(\tilde M)/T^k\to S(\tilde M)/T^k,$$
by $\psi_{i,j}(\tilde y,r)=\op{cm}(S(\tilde y,r))$,
 the center of mass of $S(\tilde y,r)$ in $\tilde y\times \Bbb R^m$. 
It is clear that $\psi_{i,j}$ is a $S(n)$-conjugate. By 
(2.7.3) and Theorem 1.3 (or a verification similar to the proof of (2.8.2) below), 
$\psi_{i,j}$ is a diffeomorphism that is close and isotopy to $\op{id}_{\tilde Y\times \Bbb R^m}$.

Consider 
$$\hat N_{i,j-1}/\hat N_{i,j}\to ((S(\tilde M)/T^k)/\hat N_{i,j},S(n)_i)\to ((S(\tilde M)/T^k)/\hat N_{i,j-1},S(n)_i).$$ 
By Lemma 2.5, $(S(\tilde M)/T^k)/\hat N_{i,j-1}\cong \tilde Y\times (\hat N_{i,0}/\hat N_{i,j-1})\cong \tilde Y\times 
\Bbb R^{m_{j-1}}$, and by induction we may assume that $(\tilde Y\times \Bbb R^{m_{j-1}},S(n)_i)$ are all $S(n)$-conjugate
by diffeomorphisms that are close and isotopy to $\op{id}_{\tilde Y\times \Bbb R^{m_{j-1}}}$. By further identifying the $S(n)_i$-actions via the conjugate diffeomorphisms,  we are back to a situation 
similar to one in the proof of $s=0$, and thus we can modify $\op{id}_{(S(\tilde M)/T^k)/\hat N_{i,j}}$ to a $S(n)$-conjugation 
via the center of mass technique. 
\qed\enddemo

\proclaim{Lemma 2.8}  Let $T^k\to (S(\tilde M),S(n)_i)@> \op{proj}_j>> (\tilde Y\times \Bbb R^m, S(n))$ be a sequence of 
$S(n)$-invariant principal $T^k$-bundles in Lemma 2.7 such that all derivative of $\tilde g_i$ have uniform bounded norm,
and all $d_{\tilde g_i}$ are close on a compact set $C$ that contains all fundamental domains of $S(M_i)$.
Then all the $S(n)_i$-actions are conjugate by diffeomorphisms that are close and isotopy to $\op{id}_{S(\tilde M)}$.
\endproclaim

\demo{Proof}  By Lemma 2.7,  from (2.6) we obtain the following $S(n)$-invariant principal $A$-bundles, 
$$\CD T^k@= T^k\\
@VV V   @VV V\\
(S(\tilde M),S(n)_i)@>\op{id}_{S(\tilde M)}>> (S(\tilde M),S(n)_j)\\
@VV \op{proj} V @VV \op{proj}V\\
(\tilde Y\times \Bbb R^m,S(n))@> \op{id}_{\tilde Y\times \Bbb R^m}>>(\tilde Y\times \Bbb R^m,S(n)).
\endCD$$
Indeed, the proof in Lemma 2.7 will go through if $\op{vol}(B_1(T^k(\tilde p),\tilde g_i))\ge v>0$ for all $\tilde g_i$ (note that $\op{vol}(B_1(T^k(\tilde p),\tilde g_i))\ge v$ implies that $\op{vol}(B_1(T^k(\tilde x),\tilde g_i))\ge c(n)^{-1}v$, see (0.3.3)). If $\op{diam}(T^k(\tilde p),\tilde g_i)\to 0$,
then this case is similar to the proof of Theorem 0.2. The remaining case to address is that $\op{vol}(B_1(T^k(\tilde p)),\tilde g_i)\to 0$ and $\op{diam}(T^k(\tilde p_i),\tilde g_i)\ge d>0$ (see Example 2.4), where it is possible that $(T^k,e,\tilde g_i)
@>\op{GH}>>(\Bbb R^t\times T^l,0)$, $t+l<k$. 

We shall apply the center of mass technique on the Riemannian universal cover of $T^k(\tilde x)$ where (2.7.3) holds. 
For any $\tilde x\in S(\tilde M)$, we define map $\nu: S(n)\to
T^k(\tilde x)$, $\phi(s)=s_j\cdot s_i^{-1}\tilde x$. Because the convexity
radius of $T^k(\tilde x)$ may be very small, we consider the
lifting map
$$\CD
S(n) @>\tilde \nu>> (\Bbb R^k,0)\\
 @|   @VV \pi V\\
S(n)@> \nu >> (T^k(\tilde x),\tilde x),
\endCD$$
where $\pi: (\Bbb R^k,0)\to (T^k(\tilde x),\tilde x)$ is the Riemannian
universal covering map and $\tilde \nu$ is the lifting map of
$\nu$ at some $\tilde y\in \pi^{-1}(\tilde x)$.
Then $\tilde \nu$
has a center of mass, denoted by $\op{cm}(\tilde x)$. We observe the following properties:

\noindent (2.8.1) For $\alpha\in \pi_1(T^k(\tilde x))$, $\op{cm}(
\alpha\cdot \tilde x)=\alpha \op{cm}(\tilde x)$ since
$\alpha$ is an isometry.

Put $\op{cm}(\tilde x)= \pi(\op{cm}(\tilde y))$. By (2.8.1),
$\op{cm}(\tilde x)$ is independent of the choice of $\tilde y$, and thus we
may formally denote it by $\op{cm}(\tilde x)=\op{cm}(s\to
s_j^{-1}s_i\cdot \tilde x))$.

\noindent (2.8.2) For $\tilde y=t(\tilde x), t\in T^k$,
$\op{cm}(\tilde y)=\op{cm}(s\to s_js_i^{-1}(t(\tilde x))) =\op{cm}(s\to
s_j^{-1}(t(s_i\tilde x)))=\op{cm}(s\to
s_j^{-1}(\phi_{i,j}(t)(s_i \tilde x)))=\op{cm}(s\to \phi_{i,j}(t)(
s_j^{-1}s_i(\tilde x)))=\phi_{i,j}(t)\cdot \op{cm}(s\to s_j^{-1}s_i(\tilde x)) =\phi_{i,j}(t)\cdot
\op{cm}(\tilde x)$, where $\phi_{i,j}: T^k\to T^k$ is the automorphism from the $T^k$-bundles weakly equivalence.

It is clear that $\hat f$ is $S(n)$-conjugate. Note that (2.8.2) implies that $\hat f:
S(\tilde M)\to S(\tilde M)$ by $\hat
f(\tilde x)=\op{cm}(\tilde x)$ is a diffeomorphism that is close and isotropy to $\op{id}_{S(\tilde M)}$.
\qed\enddemo

Putting the above work together, one obtains a proof of Theorem A. For convenience of readers, we 
summarize the proof below.

\demo{Proof of Theorem A}

We start with a collapsing sequence of $M_i$ in Theorem A,  
$M_i@>GH>>X$, and we will show that passing to a subsequence, the sequence of the Riemannian universal cover, 
 $\tilde M_i$, are all diffeomorphic. This will be done in two steps: (i) all $S(\tilde M_i)$ are diffeomorphic, (ii) all 
 $\tilde M_i$ are diffeomorphic.
 
By Theorem 0.3, we obtain the following commutative diagram,
$$\CD  F_i@> >> (S(M_i), S(n)_i)@> f_i>\op{eqGH}>(Y, S(n))\\
@ VV \op{proj}_iV    @VV \op{proj}_iV @VV \op{Proj}V\\
\bar F_i@> >>M_i@>\bar f_i >>X=Y/S(n),
\endCD$$
where $f_i$ is an $S(n)$-invariant fiber bundle map and $F_i$ is a nilpotent manifold. By Lemmas 2.2 and 2.5, 
we obtain the following $S(n)$-invariant principal connected abelian Lie group bundle, 
$$T^k\to (S(\tilde M_i),S(n)_i)@>\op{proj}_i>>(\tilde Y\times \Bbb R^m,S(n)_i).$$
Because $\pi_1(S(\tilde M_i))=\pi_2(S(\tilde M_i))=0$, by Lemma 0.4 the above principal $T^k$-bundles are all weakly equivalent; 
in particular all $S(\tilde M_i)$ are diffeomorphic.  Finally, by Lemmas 2.7 and 2.8 we conclude that all
 $(S(\tilde M_i),S(n)_i)$ are $S(n)$-conjugate, and thus all $\tilde M_i$  are diffeomorphic.
\qed\enddemo

\vskip4mm

\head 3. Appendix. The Universal Cover Diffeomorphism Finiteness of Manifolds Of Non-negative Ricci Curvature
\endhead 

\vskip4mm

In the introduction, the following problem raises: for the class of compact $n$-manifolds with uniform bounded sectional curvature and diameter, whether the Riemannian universal covers which are uniformly non-collapsed contain a finite number of 
diffeomorphic types. 

There likely is a negative answer if one relaxes the curvature to bounded Ricci curvature. Our goal here is to present 
some universal cover diffeomorphism (resp. homeomorphism) finiteness for the following class of compact $n$-manifolds of 
non-negative Ricci curvature (resp. sectional) curvature. 

\proclaim{Theorem 3.1} Given $n, v>0$, there are constants,  $\epsilon(n), C(n), C(n,v)>0$, such that if a compact $n$-manifold $M$ of $\op{diam}(M)=1$ satisfies one of the following conditions, its Riemannian universal cover $\tilde M$ contains at most $C(n)$ (resp. $C(n,v))$ many diffeomorphic (resp. homeomorphic) types respectively in (3.1.1) (resp. (3.1.2)):

\noindent {\rm (3.1.1)} $\op{Ric}_M\ge 0$ and $\op{vol}(B_{10^{-1}}(\tilde x))\ge (1-\epsilon(n))
\op{vol}(\b B_{10^{-1}}^n(0))$, $\forall\, \tilde x\in \tilde M$,  where $\b B^n_r(0)$ denotes an Euclidean $n$-ball of radius $r$.

\noindent {\rm (3.1.2)} $\op{sec}_M\ge 0$ and $\op{vol}(B_1(\tilde p))\ge v$, for some $\tilde p\in \tilde M$.
\endproclaim

Note that manifolds in (3.1.2) actually contains finitely many diffeomorphic types in the case that 
$\pi_1(M)$ contains no subgroup $\Bbb Z^{n-4}$ of finite index (see the proof of (3.1.2)).

We will review results that will be used in the proof of Theorem 3.1.

In Riemannian geometry, the celebrated splitting theorem of Cheeger-Gromoll asserts:

\proclaim{Theorem 3.2} {\rm (\cite{CG})} Let $M$ be a complete $n$-manifold of 
$\op{Ric}_M\ge 0$.

\noindent {\rm (3.2.1)} If $M$ contains a line, then $M$ splits, $M=N\times \Bbb R$.

\noindent {\rm (3.2.2)} If $M$ is compact, then the Riemannian universal cover splits, 
$\tilde M=S\times \Bbb R^k$ ($k\ge 0$), where $S$ is compact. In particular, $\pi_1(M)$
contains a free abelian group of finite index.
\endproclaim

Note that Theorem 3.2 does not provide any information on a compact manifold of $\op{Ric}_M\ge 0$ 
and finite fundamental group (equivalently, $k=0$).

Note that in a proof of Theorem 3.2, one may consider a sequence of compact $n$-manifolds satisfying the conditions 
in Theorem 3.1, and prove that for $i$ large, all $\tilde M_i$ are diffeomorphic or homeomorphic. By Theorem 3.2, 
$\tilde M_i=N_i\times \Bbb R^{k_i}$, it reduces to show that the sequence of non-collapsed compact factor $N_i$ 
are all diffeomorphic or homeomorphic (passing to a subsequence).

We recall the following well known homeomrphic and diffeomorphic stability results. 

\proclaim{Theorem 3.3} {\rm (Homeomorphic stability, \cite{Pe})} Let a sequence of compact $n$-manifolds, 
$M_i@>\op{GH}>>X$ such that 
$$\op{sec}_{M_i}\ge -1,\quad \op{diam}(M_i)\le d,\quad \op{vol}(M_i)\ge v>0.$$
Then for $i$ large, there is a homeomorphism from $M_i$ to $X$ that is also an $\epsilon_i$-GHA.
\endproclaim

Note that Theorem 3.3 was proved for a larger class that may not be topological manifolds (\cite{Pe}). 

\proclaim{Theorem 3.4} {\rm (Diffeomorphism stability, \cite{CC1})} Let a sequence of compact $n$-manifolds, 
$M_i@>\op{GH}>>X$, such that 
$$\op{Ric}_{M_i}\ge -(n-1),\quad \op{diam}(M_i)\le d,\quad \op{vol}(B_\rho(x_i))\ge (1-\epsilon(n))\op{vol}(\b B_\rho^n(0)),\,\forall\, x_i\in M_i.$$
Then $X$ is homeomorphic to a manifold $M$ and for $i$ large, there is a diffeomorphism from $M$ to $M_i$ 
such that the pullback distance functions bi-H\"ohlder converge. 
\endproclaim

We point it out that the original proof of Theorem 3.4 in \cite{CC1} uses the Reifenberg method from the geometric measure theory. Recently, it was improved that if $k=n$, then a $\delta$-splitting map is non-degenerate (\cite{CJN}). With this regularity, 
in \cite{Hu} it was showed that the differentiable map constructed in the proof of Lemma 3.7 below is 
non-degenerate, and thus is a diffeomorphism.

To apply Theorems 3.3 and 3.4, we need the following lemma.

\proclaim{Lemma 3.5} Given $n$, there exists a constant, $d(n)>0$, such that
if $M$ is a compact $n$-manifold such that
$$\op{Ric}_M\ge 0,\quad \op{diam}(M)=1,$$
and the Riemannian universal cover splits, $\tilde M=S\times \Bbb R^k$, then
$\op{diam}(S)\le d(n)$.
\endproclaim

When $k=0$ i.e., $\pi_1(M)$ is finite,  Lemma 3.5 is known (\cite{KW}). In the proof, we need the following 
generalized Margulis lemma.

\proclaim{Theorem 3.6}  {\rm (Margulis Lemma, \cite{KW})} For $n\ge 3$, there exist constants, $\epsilon(n), w(n)>0$,  such that
if a compact $n$-manifold $M$ satisfies  
$$\op{Ric}_M\ge -(n-1),\quad \op{diam}(M)<\epsilon(n),$$
then the fundamental group $\pi_1(M)$ contains a normal nilpotent subgroup of index $\le w(n)$.
\endproclaim

In the proof of Lemma 3.5, we will also use the following fact:

\proclaim{Lemma 3.7}  For $i$ large, there is a smooth map, $f_i: M_i\to N$, such that $f_i$ is also an 
$\epsilon_i$-GHA, $\epsilon_i\to 0$.
\endproclaim

To construct a smooth map, we will employ the technique of $\delta$-splitting maps in \cite{CC1}.

\definition{Definition 3.8} {\rm ($\delta$-splitting maps)} Let $M$ be a complete $n$-manifold.
For $\delta>0, r>0$, a map, $u=(u_1,...,u_k): B_r(x)\to \Bbb R^k$, is called an
$\delta$-splitting map, if
$$\Delta u_i=0,\quad |\nabla u_i|\le 1+\delta,\quad -\kern-1em\int_{B_r(x)}|g(\nabla u_i,\nabla u_j)-\delta_{i,j}|^2<\delta,\quad r^2-\kern-1em\int_{B_r(x)}|\op{Hess}_{u_i}|^2<\delta.$$
\enddefinition

\proclaim{Theorem 3.9} {\rm (\cite{CC1})} For any $\delta>0$, there is a $\epsilon(n,\delta)>0$, such that if 
$$\op{Ric}_{B_4(p)}\ge -(n-1)\epsilon,\quad d_{GH}(B_4(p), \b B^m_4(0))<\epsilon,$$
then there is an $\delta$-splitting map, $\phi=(\phi_1,...,\phi_m): B_1(p)\to \Bbb R^m$.
\endproclaim

\demo{Proof of Lemma 3.7}

Fixing any $\epsilon>0$, without loss of the generality by rescaling we may assume that the convexity 
radius of $N$ is $\ge 10$, and 
$$\op{Ric}_{M_i}\ge -(n-1)\epsilon,\quad d_{GH}(B_4(x),\b B_4^m(0))<\epsilon,\quad \forall\, x\in N,\,\,m=\dim(N),$$

Let $\{\bar x_j\}$ be an $1$-net on $N$. Then $\{B_1(\bar x_j)\}$ form a local finite cover for $N$. 
For each $j$, let $x_{i,j}\in M_i$ such that $x_{i,j}\to \bar x_j$. For $i$ large, we may assume that $\{B_2(x_{i,j})\}$ form
a local finite open cover for $M_i$.  Applying Theorem 3.9, we may assume a Harmonic map, $u_{i,j}: B_2(x_{i,j})\to \Bbb R^m$.
Then $f_{i,j}=\exp_{\bar x_j}\circ u_{i,j}: B_2(x_{i,j})\to B_2(\bar x_j)$. Take a partition of unity, $\{h_j\}$, associate to 
$\{B_2(\bar x_j)\}$. For each $x_i\in M_i$, we define an energy function on $N$:
$$E(x_i)=\frac 12\sum_jh_jd^2(f_{i,j}(x_i),y),\quad y\in N.$$
Because $E(x_i)$ is a linear combination of a finite number of strictly convex functions, $E(x_i)$ is strictly convex. 
Consequently, $E(x_i)$ achieves an unique minimum, say $y_i$. We then define $\phi_i(x_i)=y_i$. It is clear 
that $\phi_i$ is a differentiable onto map and a $\Psi(\epsilon)$-GHA.
\qed\enddemo

\demo{Proof of Lemma 3.5} 

By Theorem 3.2, the Riemannian universal cover splits, $\tilde M=S\times \Bbb R^k$ ($k\ge 0$), and 
$S$ is a compact simply connected manifold of non-negative Ricci. 

We now proceed the proof by contradiction, assuming a sequence of compact $n$-manifolds $M_i$ 
such that
$$\op{Ric}_{M_i}\ge 0,\quad \op{diam}(M_i)=1,\quad d_i=\op{diam}(S_i)\to \infty,$$
where $\tilde M_i=S_i\times \Bbb R^{k_i}$. Passing to a subsequence, we may assume 
the following commutative diagram:
$$\CD (d_i^{-1}(S_i\times \Bbb R^k),\tilde p_i, \Gamma_i)@>\op{eqGH}>> (\tilde Y\times \Bbb R^k,\tilde p,G)\\
@VV \pi_i V @VV \op{proj} V\\
(d_i^{-1}M_i,p_i)@>\op{GH}>> \op{pt}=(\tilde Y\times \Bbb R^k)/G,
\endCD$$
where $\Gamma_i=\pi_1(M_i)$.  By \cite{CoN}, $G$ is a closed Lie group of isometries, and thus $G=G_0$. 
By Theorem 3.6, $G$ is a nilpotent Lie group. Because $G$ acts transitively and effectively on $\tilde Y\times \Bbb R^k$, 
we may identify $\tilde Y$ as a compact nilpotent subgroup, and thus $\tilde Y=T^s$, a torus. 
Hence, $d_i^{-1}S_i@>\op{GH}>>T^s$ with $\op{Ric}_{S_i}\ge 0$ and $\op{diam}(T^s)=1$. 

By Lemma 3.7, for $i$ large there is a smooth onto map, $f_i: d_i^{-1}S_i\to T^s$, which is also an $\epsilon_i$-GHA, $\epsilon_i\to 0$. For any closed geodesic $\alpha$ in $T^s$, $\alpha$ is not homotopically trivial, then there is a loop $\gamma_i$  
in $S_i$ such that $f_i(\gamma_i)$ is homotopically equivalent to $\alpha$; a contradiction because $S_i$ is 
simply connected. 
\qed\enddemo

\demo{Proof of Theorem 3.1} 

By Theorem 3.2, $\tilde M$ splits, $\tilde M=S\times \Bbb R^k$, $0\le k\le n$, and by Lemma 3.5, 
$\op{diam}(S)\le d(n)$. It suffices to consider a sequence of compact $m$-manifolds, $S_i@>\op{GH}>>
X$ such that $\op{diam}(S_i)\le d(n)$ and  $\op{Ric}_{S_i}\ge 0$ or $\op{sec}_{S_i}\ge 0$ respectively, and 
show that for $i$ large, all $S_i$ are all diffeomorphic or homeomorphic respectively.

(3.1.1) The volume condition is equivalent to an uniform Reifenberg condition, and by Theorem 3.4 
for $i$ large $S_i$ is diffeomorphic to $S_{i+j}$ for all $j\ge 1$. 

(3.1.2) Note that $\op{vol}(B_1(\tilde p))\ge v$ implies that $\op{vol}(B_1(S,\tilde p))\ge \epsilon(n,v)>0$. 
By Theorem 3.3, for $i$ large all $S_i$ are homeomorphic (the homeomorphic finiteness was first obtained 
by \cite{GPW1,2}). 
For $\dim(S)\ne 4$,  the homeomorphic stability implies a diffeomorphism finiteness (\cite{KiS}). 
For $\dim(S)=4$ and thus $n\ge 4$, $\pi_1(M)$ contains a subgroup $\Bbb Z^{n-4}$ of finite index.
\qed\enddemo


\vskip10mm

\Refs
\nofrills{References}
\widestnumber\key{APS12}

\vskip3mm

\ref \key An \by M. Anderson\pages 429-445 \paper Convergence and
rigidity of manifolds under Ricci curvature bounds\jour Invent.
math. \vol 102 (2) \yr 1990
\endref

\ref \key AC \by M. Anderson; J. Cheeger
\pages  265-281
\paper $C^\alpha$-compactness for manifolds with Ricci curvature and injectivity radius  bounded below
\jour J. Diff. Geom
\vol 35
\yr 1992
\endref

\ref
\key Ch
\by J. Cheeger
\pages 61-75
\paper Finiteness theorems for Riemannian manifolds
\jour Amer. J. Math
\vol 92
\yr 1970
\endref

\ref
\key CC1
\by J. Cheeger; T. Colding
\pages 189-237
\paper Lower Bounds on Ricci Curvature and the Almost Rigidity of Warped Products
\jour Ann. of Math.
\vol 144
\issue 1
\yr Jul., 1996
\endref

\ref
\key CC2
\by J. Cheeger; T. Colding
\pages 406-480
\paper On the structure of space with Ricci curvature bounded below I
\jour J. Diff. Geom
\vol 46
\yr 1997
\endref

\ref
\key CFG
\by J. Cheeger, K. Fukaya; M. Gromov
\pages 327-372
\paper Nilpotent structures and invariant metrics on collapsed manifolds
\jour J. Amer. Math. Soc.
\vol 5
\yr 1992
\endref

\ref
\key CG1
\by J. Cheeger; M. Gromov
\pages 309-346
\paper Collapsing Riemannian mamfolds while keeping their curvature bounded I
\jour J. Differential Geom.
\vol 23
\yr 1986
\endref

\ref
\key CG2
\by J. Cheeger; M. Gromov
\pages 269-298
\paper Collapsing Riemannian mamfolds while keeping their curvature bounded II
\jour J. Differential Geom.
\vol 32
\yr 1990
\endref

\ref \key CG \by J. Cheeger; D. Gromoll \pages 119-128\paper
The splitting theorem for manifolds of nonnegative Ricci curvature
 \jour J. Diff. Geom. \yr 1971 \vol 6
\endref

\ref
\key CJN
\by J. Cheeger,  W. Jiang; A. Naber
\pages
\paper Rectifiability of singular sets in noncollapsed spaces with Ricci curvature bounded below
\jour arXiv:1805.07988
\vol
\yr 
\endref

\ref
\key CoN
\by T. Colding; A. Naber
\pages 1172-1229
\paper Sharp H\"older continuity of tangent cones for spaces with a lower Ricci curvature bound and applications
\jour Ann. of Math
\vol 176
\yr 2012
\endref

\ref
\key Es
\by J.-H Eschenburg
\pages 469-480
\paper New examples of manifolds with strictly positive curvature
\jour Invent. Math
\vol 66 
\yr 1982
\endref

\ref \key FR1 \by  F. Fang; X. Rong \pages 641-674 \paper Positive
pinching, volume and homotopy groups \jour Geom. Funct. Anal \yr
1999 \vol 9
\endref

\ref \key FR2 \by  F. Fang; X. Rong \pages 61-109 \paper The
twisted second Betti number and convergence of collapsing
Riemannian manifolds \jour Invent. Math. \yr 2002 \vol 150
\endref

\ref
\key Fu1
\by K. Fukaya
\pages139-156
\paper Collapsing of Riemannian manifolds to ones of lower dimensions
\jour J. Diff. Geom.
\vol25
\yr1987
\endref

\ref
\key Fu2
\by K. Fukaya
\pages1-21
\paper A boundary of the set of Riemannian manifolds with bounded curvature and diameter
\jour J. Diff. Geom.
\vol28
\yr1988
\endref

\ref
\key Fu3
\by K. Fukaya
\pages333-356
\paper Collapsing of Riemannian manifolds to ones of lower dimensions II
\jour Jpn. Math.
\vol41
\yr1989
\endref

\ref
\key FY
\by K. Fukaya; T. Yamaguchi
\page 253-333
\paper The fundamental groups of almost non-negatively
curved manifolds
\jour Ann. of Math
\vol 136
\yr 1992
\endref

\ref
\key Hu
\by H. Huang
\pages
\paper Fibrations, and stability of compact Lie group actions on manifolds with local bounded Ricci covering geometry
\jour  Preprint
\vol
\yr
\endref

\ref \key GW \by R. E. Green; H. Wu\pages 119-141\paper Lipschitz
convergence of Riemannian manifolds\jour Pacific J. Math. \vol
131\yr 1988
\endref

\ref\key Gr\by M. Gromov\pages 231-241\paper Almost flat manifolds\jour J. Diff. Geom.\vol 13\yr 1978\endref

\ref \key GLP \by M. Gromov, J. Lafontaine; P. Pansu \paper
Structures metriques pour les varietes riemannienes \jour
CedicFernand Paris, \yr 1981
\endref

\ref \key GK\by K. Grove; H. Karcher\pages 11-20\paper How to
conjugate $C^1$-close group actions \jour Math. Z. \vol 132\yr
1973
\endref

\ref
\key GPW1
\by K. Grove; P. Petersen; J. Wu
\pages 205-213
\paper Controlled topology in geometry
\jour Invent. Math
\vol 99
\yr 1990
\endref

\ref
\key GPW2
\by K. Grove; P. Petersen; J. Wu
\pages 221-222
\paper Controlled topology in geometry (Erraturm)
\jour Invent. Math
\vol 104
\yr 1991
\endref

\ref
\key KiS
\by P. Kirby; L. Siebenmann
\pages
\paper Foundational essays on topological manifolds, smoothings and triangulations
\jour Ann. Math. Stud.,  Princeton University Press
\vol 88
\yr 1977
\endref

\ref
\key KoS
\by O. Kowalsk; M. Sekizawa
\pages 1799-1809
\paper On the geometry of orthonormal frame bundles
\jour Math. Nachr.
\vol 281 
\yr 2008
\endref

\ref
\key KW
\by V. Kapovitch; B. Wilking
\pages
\paper Structure of fundamental groups of manifolds with Ricci curvature bounded below
\jour Preprint
\vol 
\yr 2011
\endref

\ref
\key MT
\by M. Mimura; H. Toda
\pages
\paper  Topology of Lie groups I and II
\jour Amer. Math. Socie. Providence, Rhode Island
\vol 91
\yr 1991
\endref

\ref
\key Pa
\by R.S. Palais 
\pages 362-364
\paper Equivalence of nearby differentiable actions of a compact group
\jour  Bull. Amer. math. Soc.
\vol 67
\yr 1961
\endref

\ref \key Pe\by G. Perelman \pages
\paper Alexandrov spaces with curvatures bounded from below II\jour
preprint \vol \yr
\endref

\ref \key Pet \by S. Peters \pages 77-82\paper Cheeger's
finiteness theorem for diffeomorphism classes of Riemannian
manifolds\jour J. Reine Angew. Math \vol 349 \yr 1984
\endref

\ref \key PT\by A. Petrunin; W. Tuschmann\pages 736-774\paper
Diffeomorphism finiteness, positive pinching, and second homotopy
\jour Geom. Funct. Anal. \vol 9\yr 1999
\endref

\ref
\key Ro
\by X. Rong
\pages 193-298
\paper Convergence and collapsing theorems in
Riemannian geometry
\jour Handbook of Geometric Analysis, Higher
Education Press and International Press, Beijing-Boston
\issue II ALM 13
\yr 2010
\endref

\endRefs
\enddocument